\documentclass[12pt,oneside]{article}
\usepackage{amsmath,amsthm,amsfonts,amssymb,MnSymbol,mathrsfs,latexsym}
\usepackage{mathtools}

 \newcommand{\const}{\rm const}

\theoremstyle{plain}

\newtheorem{remark}{Remark}[section]

\renewenvironment{proof}{{\bf{Proof.}}}{\hfill $\Box$ \\}

\title{\large \textbf{Confidence regions for the multidimensional density  in the uniform norm based on the recursive
Wolverton-Wagner estimation}}

\footnotesize\date{}

\author{\normalsize Maria Rosaria Formica ${}^{1}$,   \normalsize Eugeny Ostrovsky
${}^2$ and \normalsize Leonid Sirota ${}^3$}

\begin{document}

  \maketitle

\begin{center}
{\footnotesize ${}^{1}$ Universit\`{a} degli Studi di Napoli \lq\lq Parthenope\rq\rq, via Generale Parisi 13,\\
Palazzo Pacanowsky, 80132,
Napoli, Italy.} \\

\vspace{1mm}

{\footnotesize e-mail: mara.formica@uniparthenope.it} \\

\vspace{2mm}

{\footnotesize ${}^{2,\, 3}$  Bar-Ilan University, Department of Mathematics and Statistics, \\
52900, Ramat Gan, Israel.} \\

\vspace{1mm}

{\footnotesize e-mail: eugostrovsky@list.ru}\\

\vspace{1mm}

{\footnotesize e-mail: sirota3@bezeqint.net} \\

\end{center}

\vspace{4mm}

\begin{abstract}
\hspace{3mm}  We construct an optimal {\it exponential tail decreasing} confidence region for an unknown density of distribution in the
Lebesgue-Riesz as well as in the {\it uniform} norm, built on the sample of the  random vectors based of the famous {\it recursive}
Wolverton-Wagner density estimation.
\end{abstract}

\vspace{4mm}

 {\it \footnotesize Keywords:}

 \vspace{3mm}

 {\footnotesize  Probability space, random variables and processes (fields), uniform  and Lebesgue - Riesz norms and spaces, density, bias,
Wolverton - Wagner and Parzen - Rosenblatt estimations, Rosenthal's inequality, sample, independence, ordinary Euclidean space, loss function,
 Chernoff's  exponential tail  estimate, Lebesgue - Riesz norm and spaces, smooth functions, distribution and tail of distribution,
confidence regions,  Young - Fenchel (Legendre) transform, conditions of orthogonality, bias and variation, convex domain,  Lipschitz condition,
moments, kernel, windows or bandwidth, multivariate Banach spaces of random variables,  exponential estimates, Grand Lebesgue Spaces.}

\vspace{5mm}

 \section{Introduction}

\vspace{4mm}

 \hspace{3mm} Let $ \ (\Omega, \cal{F}, {\bf P}) \ $ be a non-trivial probability space with an expectation $ \ {\bf E} \ $ and a variance $ \ {\bf Var}. \ $
Let also $ \ \xi_1, \xi_2, \ldots,\xi_k, \ldots, \ \xi_n; \ k \in \mathbb N = 1,2,\ldots,n; \ n \ge 2 \ $ be a {\it sample}, i.e. finite sequence of independent,
identical distributed, \ (i, i.d.) \ random vectors (r.v.) taking the values in the ordinary Euclidean space
$ \ \mathbb R^d, d \geq 1 \ $, equipped with usually distance function $ \  ||x - y|| = (x - y, x - y)^{1/2}, \ x,y \in R^d, \ $
and having  an unknown distribution {\it density}  function
$ \ f = f(x), \ x \in \mathbb R^d. \ $  C. Wolverton and T.J. Wagner in \cite{Wolverton Wagner}  introduced the famous statistical estimation
$ \ f_n^{ W W}(x) = f_n(x) \ $ for the function $ \ f(\cdot)$ as follows. \par

 \hspace{3mm} Let $ \{h_k\}, \ k\in \mathbb N $, be a positive deterministic sequence of real numbers such that
$ \lim_{k \to \infty} h_k = 0, \ h_1 = 1, \ $ ({\it windows, or bandwidth}). Let also  $ K = K(x)$,  $x \in \mathbb R^d $, be
 certain {\it kernel}, i.e. a measurable numerical valued fixed {\it "sufficiently smooth,"} see an exact definition further,
 even normalized function, i.e. for which

$$
\int_{\mathbb R^d} K(x) dx = 1.
$$

\vspace{3mm}

\hspace{3mm} Recall that following definitions.

\vspace{3mm}

\noindent {\bf  Definition of the Wolverton - Wagner recursive estimate}, see \cite{Wolverton Wagner}; as well as
\cite{Nadaraya_Babilua}.

\vspace{3mm}

\begin{equation} \label{def WW}
f_n^{ W W}(x) = f_n(x) \stackrel{def}{ =} \frac{1}{n} \sum_{k=1}^n \frac{1}{h_k^d} K \left( \ \frac{ x - \xi_k}{h_k}   \right).
\end{equation}

\vspace{3mm}

\noindent {\bf  Definition of the Parzen - Rosenblatt (or Kernel) estimate} (see, e.g., \cite{Parzen,Rosenblatt}).

\begin{equation} \label{def PR}
g_n^{P R}(x) = g_n(x)\stackrel{def}{ =} \frac{1}{n \ h_n^d} \sum_{k=1}^n K \left( \ \frac{ x - \xi_k}{h_n} \right).
\end{equation}

\vspace{3mm}

\noindent Note that the Wolverton - Wagner estimate obeys a very important recursion property:

$$
f_{n + 1}^{WW}(x) = \frac{n}{n + 1} f_n^{WW}(x) + \frac{h^{-d}_{n + 1}}{n + 1} \cdot K \left(\ \frac{x - \xi_{n + 1}}{h_{n + 1}} \ \right).
$$

\vspace{3mm}

 \ The recurrent definition of probability density estimates $ \ f_n^{W W}(x) \ $ has two
obvious advantages: 1) there is no need to memorize data, i.e. if the estimate
$ \ f_n^{W W} (x) \ $  is known, then the following one $ \ f_{n+1}^{W W}(x)  \ $ can be calculated by means of the last observation
$ \ f_n^{WW} \ $  only, without using the sampling $ \ \xi_1, \xi_2,\ldots, \xi_n$; 2) the asymptotic variation as well as {\it bias}
of the estimate $ \ f_n^{W W}(x) \ $ does not exceed the asymptotic variation and bias  of the estimate $ \ g_n^{P R}(x), \ $ see e.g.
\cite{Devroye 1,Devroye 2,Form Os Sir3}. This proposition holds true still when the error of estimation
is understood in the classical Lebesgue - Riesz norm sense

$$
R_p[f, f_n] \stackrel{def}{=} {\bf E} \left\{ \int_{\mathbb R^d} |f^{WW}_n(x) - f(x)|^p \ dx \ \right\}^{1/p}, \ p \in [1,\infty).
$$

\vspace{3mm}

 \ The {\it loss function} is understood in this report in an {\it uniform norm} deviation

\vspace{3mm}

\begin{equation} \label{unif norm}
R_{\Psi, B_n, D}[f, f_n] \stackrel{def}{=} {\bf E} \Psi \left[ \ B_n \ \sup_{x \in D} |f^{WW}_n(x) - f(x)| \ \right],
\end{equation}

\vspace{3mm}

where $ \ \Psi(\cdot) \ $ is certain {\it weight} function; indeed  some   non - negative  Young function. \par
 \ This means by definition that this function is such that $ \ \Psi: [0,\infty) \to [0,\infty),  \ $ is {\it strictly} increasing,
continuous, convex, and

$$
\Psi(0) = 0; \hspace{3mm} \lim_{x \to \infty} \Psi(x) = \infty.
$$

\vspace{3mm}

 \  {\sc Here and hereafter}  $ \  D \ $ {\sc is fixed non - empty closed convex bounded subdomain (compact) subset of the whole space} $ \ R^d. \ $ \par

\vspace{3mm}

 \ $ \ B_n \ $ is some deterministic normed positive numerical sequence such that $ \ \lim_{n \to \infty} B_n = \infty. \ $
 \ For instance,  $ \ \Psi(y) = |y|^m, \ m = \const > 1 \ $ or $ \ \Psi(y) = \exp(|y|^m) - 1, \ m = \const > 0. \ $ \par

\vspace{3mm}

\hspace{3mm}  {\bf Our aim in this report is to deduce the exact exponential decreasing
estimate for the tail of deviation probability} \ $ \ S^{WW}(u) = S^{WW}_{D,n,d}(u) \stackrel{def}{=} \ $

\vspace{3mm}

\begin{equation} \label{aim}
 {\bf P}(B_n \sup_{x \in D} |f^{WW}_n(x) - f(x)| > u); \ u > u_0 = \const > 0,
\end{equation}

\vspace{3mm}
of course under  appropriate restrictions  and for exact {\it optimal} deterministic positive  normed numerical sequence $ \ B_n$,  tending to infinity.  \par
 \ To be more precise, one can suppose that the function $ \ f(\cdot) \ $ belongs to one or another smoothness class of functions $ \ U = \{f\}. \ $
Of course, the {\it optimal} norming sequence $ \ B_n \ $ depends on $ \ U: \ B_n = B_n(U). \ $   \par

\vspace{3mm}

 \ For the classical Parzen - Rosenblatt estimate  $ \ f_n^{PR}(x) \ $ analogous results were obtained, e.g., in \cite{Form Os Sir2,Form Os Sir3},
 \cite[Chapter 5, section 5.2]{Ostrovsky1999}. \par

\vspace{4mm}

 \section{Auxiliary  considerations.}

\vspace{3mm}

 \hspace{3mm} We must introduce several  notations and definitions. Let $ \ \beta \ $  be a fixed positive constant: $ \ \beta = \const \in (0,\infty). \ $
 Denote by $ \ [\beta] \ $ its integer part and correspondingly by $ \ \{\beta\} \ $ its  fraction part $ \ \{ \beta \} = \beta - [\beta]. \ $ Define as ordinary
 the following (very popular) functional class $ \ \Sigma(\beta,L)  \ $   consisting on all the continuous bounded functions such that all its partial derivatives of (integer vector)
 order $ \  \vec{\alpha} = \{\alpha_1,\alpha_2, \ldots, \alpha_d \}, \ $ where

\vspace{3mm}

 \begin{equation} \label{condition alpha 1}
 |\alpha| \stackrel{def}{=} \sum_{j=1}^d  \alpha_j \le [\beta],
\end{equation}

\vspace{3mm}

are bounded and satisfy the H\"older  condition  with degree $ \ \{\beta\}: \ \exists L = \const < \infty $   such that, for all the values $ \ \ \alpha \ $ from the set
 (\ref{condition alpha 1}),

\vspace{3mm}

\begin{equation} \label{condition alpha 2}
\left| \frac{D^{|\alpha|} f(x)}{\prod_{j=1}^d \partial^{\alpha_j} x_j}  -  \frac{D^{|\alpha|} f(y)}{\prod_{j=1}^d \partial^{\alpha_j} y_j} \ \right|
 \le   L \ ||x - y||^{\{\beta\}}.
\end{equation}

\vspace{3mm}

 \ In the case when $ \ \beta \ $ is integer: \ $ \ \{\beta\} = 0, \ $ in the right hand of (\ref{condition alpha 2}) must stay
$ \ L \ ||x - y||;\ $ the so - called Lipschitz condition. \par
\ Set also $ \ \Sigma(\beta) := \cup_{L > 0} \Sigma(\beta,L). \  $  \ Obviously, the set $ \ \Sigma(\beta) \ $ forms a Banach  functional space. \par

\vspace{3mm}

 \hspace{3mm} {\sc  We suppose further that the unknown density function belongs to the set} $ \ \Sigma(\beta), \ $ or equally \par

\begin{equation} \label{belonging to class}
\exists \beta \ge 0, \ \exists L \in (0,\infty) \ : \ f(\cdot) \in \Sigma(\beta, L).
\end{equation}

\vspace{3mm}

 \ See the following works devoted to the consistent measurement (estimation) of the parameters $ \ (\beta,L) \ $
\cite{Goldenshluger Lepski,Khardani Slaoui,Bertin 1,Bertin 2,Devroye 1} etc.\par

\vspace{3mm}

 \ {\it Let us now impose also several conditions on the (measurable)  kernel - function}  $ \ K(\cdot). \ $

 \vspace{3mm}

$$
K(x) = K(-x);  \ \int_{\mathbb R^d} K(x) \ dx = 1; \ \int_{\mathbb R^d} K^2(x) dx < \infty, \  K(\cdot) \in L_1(\mathbb R^d) \cap C(\mathbb R^d);
$$

$$
\ \exists \delta  = {\const} \in (0,1], \  \exists c = {\const} (0, \infty) \ :  \ |K(x) - K(y)| \le c |x - y|^{\delta}.
$$

\vspace{3mm}

 \ The following imposed by us important conditions (of orthogonality) are needed only by the bias estimation: for arbitrary non
 negative integer vector $ \ l = \vec{l} \in \mathbb R^d \ $ such that  $ \ |l| \le [\beta]  \ $ we suppose that

\vspace{3mm}

\begin{equation} \label{cond bias}
\int_{\mathbb R^d} x^l \ K(x) \ dx = 0; \hspace{3mm} x^l \stackrel{def}{=} \prod_{j=1}^d x_j^{l_j}.
\end{equation}

\vspace{3mm}

$$
|x|:=  \sqrt{ \sum_{j=1}^d x_j^2}.
$$

\vspace{3mm}

 \ The optimal  chooses (in different senses) of the kernel $ \ K(\cdot) \  $ is investigated in the report  \cite{Form Os Sir3}. \par

\vspace{3mm}

 So, we allow consider in general case as a capacity of this kernels alternating ones. \par

 \ Let us discuss about the choice of the windows (bandwidth) $ \ \{h_k\}. \ $  For the classical Parzen - Rosenblatt estimate
  $ \ h_k \ $ depends only on $ \ n: $

 $$
h_n \stackrel{def}{=} c \left\{ \ \frac{\ln n}{n} \ \right \}^{\beta/(2 \beta + d)}, \ k = 2,3. \ldots,n,
 $$
 the  asymptotical optimal  choice.\par

 \vspace{3mm}

\ Let us return to the Wolverton - Wagner estimate.\par

\vspace{3mm}

\  Pick  then for the definiteness $ \ h_1 = 1. \ $ The asymptotically  as $ \ n \to \infty \ $ optimal in
 the uniform norm $ \  \sup_x |f_n(x) - f(x)| = ||f_n - f||C  \ $
   the Wolverton - Wagner  estimates $ \ f_n(\cdot) = f_n^{WW}(\cdot) \ $ of the density $ \ f(\cdot) \ $
 values of the windows  (bandwidth)  $ \ \{h_k\} \ $  under formulated above restrictions have the   following form.

 \vspace{3mm}

\begin{equation} \label{optimal wind}
 h_k \stackrel{def}{=} c_1 \left\{ \ \frac{\ln k}{k} \ \right \}^{\beta/(2 \beta + d)}, \ k = 2,3. \ldots,n;
\end{equation}

 \vspace{3mm}

see e.g.\cite{Bertin 1},  \cite{Bertin 2}, \cite{Goldenshluger  Lepski} etc. \par

\vspace{3mm}

\  {\it So, we choose furthermore the  values} $ \ \{h_k\} \ $ {\it in accordance with the relation} (\ref{optimal wind}) {\it for
the Wolverton - Wagner estimate}.\par

\vspace{3mm}

 \ Introduce the following variables $ \ B_1 = 1, \ $

$$
B_n = B_n(\beta,d) \stackrel{def}{=} \left[ \   \frac{n}{\ln n} \ \right]^{\beta /(2 \beta + d)} \asymp 1/h_n, \  n \ge 2;
$$

$$
Q_n (u)= Q_n(d; \beta,L; u) \stackrel{def}{=} {\bf P}\left( \ B_n \cdot \sup_{x \in D} |f^{WW}_n(x) - f(x)| > u \ \right), \ u \ge e.
$$

\vspace{3mm}

 \hspace{3mm}  Notice that  under formulated above  notations and conditions

$$
\sup_n \sup_{x \in D}  \  B_n(\beta,d) | \ {\bf E} f^{WW}_n(x) - f(x)| < \infty,
$$
{\it therefore it is sufficient to investigate only the variable }

$$
Q_n^o (u)= Q^o_n(d; \beta,L; u) \stackrel{def}{=} {\bf P}\left( \ B_n \cdot \sup_{x \in D} |f^{WW}_n(x) - {\bf E}f^{WW}_n(x)| > u \ \right), \ u \ge e.
$$

\vspace{3mm}

\section{Main result.}

\vspace{3mm}

 \hspace{3mm} Let $ \ (Z, \cal{Z}, \nu) \ $ be arbitrary measurable space equipped  with a sigma - finite measure $ \ \nu. \ $ Denote as ordinary the classical
 Lebesgue - Riesz space $ \  L_{p,Z,\nu} = L_p(Z), \ p \in[1,\infty)  \ $ as a set of all the measurable numerical valued functions $ \ h: Z \to R \ $
 having a finite norm

$$
||h||_{p; Z,\nu} \stackrel{def}{=} \left[ \ \int_Z |h(z)|^p \ \nu(dz)  \ \right]^{1/p}.
$$

\vspace{3mm}

 \ In particular, let $ \zeta: \Omega \to \mathbb R $ be an arbitrary measurable function (random variable), denote as usual the corresponding Lebesgue - Riesz
 $ \ L_p  = L_p(\Omega) \ $ norm

$$
||\zeta||_p \stackrel{def}{=} \left[ \ {\bf E} |\zeta|^p \ \right]^{1/p}, \ \ p \in [2,\infty).
$$

\vspace{3mm}

 \hspace{3mm} We conclude under our notations and restrictions

 \vspace{3mm}

\begin{equation} \label{Leb Riesz popint}
 \sup_{x \in D} \ || \ |f_n^{WW}(x) - f(x)| \ ||_p \le C_1(\beta,L,d,D) \times \frac{p}{\ln p} \times
\end{equation}

\vspace{3mm}

$$
\left( \ \frac{\ln n}{n} \ \right)^{\beta/(2 \beta + d)}, \ n \ge 2; \hspace{3mm}    C_1 \stackrel{def}{=} C_1(\beta,L,d) \in (0,\infty).
$$

\vspace{3mm}

\begin{proof}

 \ It follows immediately, quite similarly as in \cite{Form Os Sir2}, from the classical Rosenthal inequality for the moment for sums of independent
 centered r.v. (\cite{Rosenthal}); see also \cite{Ibragimov_Hasminskii,Ibragimov Sharachmedov 1,Ibragimov Sharachmedov 2}. \\
 \hspace{3mm} The exact value of the correspondent constant is  derived in \cite{Naimark Ostrow}.\\
\end{proof}

\vspace{2mm}

\begin{remark}
{\rm The  relation (\ref{Leb Riesz popint}) may be rewritten on the language of the so - called Grand Lebesgue Spaces (GLS), see e.g.
\cite{Ermakov etc. 1986,Kozachenko-Ostrovsky 1985,Kozachenko at all 2018,liflyandostrovskysirotaturkish2010,Ostrovsky1999,Ostrovsky HIAT,Ostrov Sir2}.}
\end{remark}

\vspace{3mm}

 \ Namely, let $ Y(\cdot)$ be an arbitrary measurable numerical valued function, for instance a random variable and, for $1 < a < b \le \infty$, let
 $ \psi(p), \ p \in (a,b)$, be a strictly positive numerical valued continuous function.  The norm in the Grand Lebesgue spaces $G\psi(a,b)$ is defined by

\vspace{3mm}

\begin{equation} \label{Norm GLS}
||Y||_{G\psi(a,b)} \stackrel{def}{=} \sup_{p \in (a,b)} \ \left\{ \frac{||Y||_p}{\psi(p)} \ \right\}.
\end{equation}

 \vspace{3mm}

 \ Let for instance $ \xi \ $ be a r.v. such that $ \ ||\xi||G\psi \le \kappa, \ 0 < \kappa < \infty. \ $ Define the auxiliary function
 $ \  h(p) = h[\psi](p) := p\ln \psi(p). \ $ and  introduce its Young - Fenchel, or Legendre transform

$$
h^*(t) = h^*[\psi](t) \stackrel{def}{=} \sup_{p \ge 2}(p t - h[\psi](p)), \ v \ge 1.
$$
 \ Another name: conjugate function. \par
 \ It is well known the following Chentzov's inequality

$$
T[\xi] (t) \le \exp \left( \  - h^*(t/\kappa) \ \right), \ t \ge \kappa.
$$
 \vspace{3mm}

 \ If for example $ \ \psi(p) = \psi_j(p) \stackrel{def}{=} \frac{p}{\ln p}, \ p \ge 2,  $ then

 \vspace{3mm}

 \begin{equation} \label{conj funk}
 h^*[\psi_l](t)  \sim  \ t  +  \ln t \cdot \ln \ln t, \ t \ge e^e.
 \end{equation}

 \vspace{4mm}

 \ Define the function
$$
\psi_l(p) \stackrel{def}{=} \frac{p}{\ln p}, \ \ p \in (2,\infty);
$$
then
\begin{equation} \label{psil}
 B_n(\beta,d) \ sup_{x D} \ || f_n^{WW}(x) - f(x) \ ||_{G{\psi_l}(2,\infty)} \le C_1(\beta,L,d,D).
\end{equation}

\vspace{3mm}

 \ This proposition follows from the last relation (\ref{psil}) a corresponding tail estimation, see e.g. \cite{Ermakov etc. 1986,Ostrovsky1999}.
 Define, as ordinarily, for an arbitrary numerical valued random variable (r.v.) $ \ \eta, \ $ its tail function
$$
T_{\eta}(t) \stackrel{def}{=} {\bf P} (|\eta| \ge t), \hspace{3mm}  t \ge t_0 = {\const} > 0.
$$
 \ Introduce  also the following family of random fields
$$
\zeta_{n}(x) :=   B_n(\beta,d) \cdot  \ \left\{ \ f_n^{WW}(x) - f(x) \ \right\}; \ \ z := t/e, \ \ z \ge 1;
$$

$$
\nu(z) \stackrel{def}{=} \exp(-z) \cdot \exp(-\ln z \cdot \ln \ln z);
$$
then

\begin{equation} \label{one dim tail}
\sup_n \ \sup_x T_{\zeta_{n,x}}(t) \le C_1 \nu(z) = C_1 \nu(t/e).
\end{equation}

\vspace{3mm}

\hspace{3mm}  Further, we deduce analogously as in (\ref{psil}),

\begin{equation} \label{Delta psil}
\begin{split}
 & B_n(\beta,d) \cdot \  ||(f_n^{WW}(x) - f(x)) - (f_n^{WW}(y) - f(y))||_{G{\psi_l}(2,\infty)}\\[2mm]
& \le C_2(\beta,L,d) \times ||x - y||/h_n, \ \ x,y \in \mathbb R^d.
\end{split}
\end{equation}

\vspace{3mm}

\ It follows from the  well - known theory of Grand Lebesgue Spaces (GLS), see e.g.  \cite{Ermakov etc. 1986},
\cite{Ostrovsky1999}, chapter 5, section 5.4 - 5.5, pages 253 \ - \ 264, the following result. \par

\vspace{3mm}

\hspace{3mm} {\bf Theorem 3.1.} \par

 \ Let $ D $ be a non - empty compact subset on the whole domain $ \ \mathbb R^d \ $. We  propose under all the formulated above
 conditions and notations

\vspace{3mm}

\begin{equation} \label{main res}
 {\bf P}\left( \ B_n(\beta,d) \cdot \sup_{x \in D} |f^{WW}_n(x) - f(x)| > t \ \right) \le  C_{12}(D,\beta,d,L) \ \nu(t/e), \ \ t \ge e,
\end{equation}

\vspace{3mm}
 where as ordinary $ \ C_{12}(D,\beta,d,L) < \infty. \ $ \par

\vspace{3mm}
 \ The {\bf proof}  is completely alike to the one for the classical Parzen - Rosenblatt estimate, see \cite{Ostrovsky1999}, chapter 5, section 5.4.
 It used a famous method belonging to  R.Reiss \ \cite{Reiss} \ of partition of the set $ \ D \ $ on the subsets of the variable  small volumes.\par

 \ In detail, introduce an auxiliary  numerical valued {\it centered} random process (field), (r.f.),   more precisely,  a family ones

$$
\zeta_n(x) \stackrel{def}{=} B_n(\beta,d) \times \left[ \ f^{WW}_n(x) - {\bf E}f^{WW}_n(x) \ \right], \ n = 1,2,\ldots; \ x \in D.
$$

  \ We must estimate the tail of maximum distribution for this r.f.

$$
{\bf P_z}(u) \stackrel{def}{=} {\bf P} \left( \ \sup_{x \in D} |\zeta_n(x)| > u \ \right),
$$
of course, for all sufficiently greatest values $ \ u, \ $ say, for $ \ u \ge e. \ $ \par
 \ One can assume without loss of generality that the set $  \ D \ $  is unit "cube"  $ \ D = [0,1)^d. \ $
 \ Put now  $ \ q = (2 \beta + d)/(\beta + d).  \ $   Introduce also the following system of cubes $ \ \vec{k} = \{ \ k_1, k_2, \ldots, k_d \ \}; \ $

 $$
 L(\vec{k}) := \otimes_{j=1}^d \left\{ \  \left[\frac{k_j}{n}, \ \frac{k_j+1}{n} \ \right)  \ \right\},
 $$
 where $ \ j = 1,2,\ldots,d; \hspace{2mm} k_j = 0, 1, 2, \ldots, n-1; \ $ which covers the whole set $ \ D. \ $ \par

 \ Define also the functions

 $$
{\bf P}_a(u) \stackrel{def}{=} \max_{\vec{k}}  {\bf P} \left[\sup_{x \in L(\vec{k})} B_n |f_n(x) - f(x)| > u \ \right].
 $$

\vspace{3mm}

\hspace{3mm}  Further, we deduce quite analogously as in (\ref{psil})

$$
 B_n(\beta,d) \cdot \  ||[f_n^{WW}(x) - {\bf E}f_n^{WW}(x)] - [f_n^{WW}(y) - {\bf E} f_n^{WW}(y)]||_{G{\psi_l}(2,\infty)} \le
$$

\vspace{3mm}

\begin{equation} \label{Delta distance 2}
 C_4(\beta,L,d) \times ||x - y||/h_n, \ \ x,y \in \mathbb R^d.
\end{equation}

\vspace{3mm}

\ Notice that the metric entropy $ \  H(D, \rho, \epsilon)  \ $
of the set $ \ D \ $ relative the distance  function $ \  \rho(x,y) := ||x - y||/h_n  \ $  may be estimated as follows

$$
 H(D, \rho, \epsilon)  \le C_1 +  C_2(d) \ln n + C_3 |\ln \epsilon|, \ \epsilon \in (0,1).
$$

 \vspace{3mm}

 \hspace{3mm}  Note  also that  it is sufficient for this purpose to ground the estimate of the form (\ref{main res})
 function for the alike   but for the  following {\it centered} correspondent random variables

\vspace{3mm}

\ Consider for definiteness the first cube $ \ Q  \stackrel{def}{=} [0, 1/n]^d,  \  $  i.e. when $ \ k_j = 0,1;  \ldots, \ j = 1,2, \ldots,d. \ $
 \ It follows immediately from the Rosenthal's moments estimates for the sums of the independent centered r.v. - s,  \ see e.g.
\cite{Ibragimov Sharachmedov 1}, \cite{Naimark Ostrow}, \cite{Rosenthal}, \  that for $ \ x \in Q \ $  and  $ \ u \ge e \ $

\vspace{3mm}

\begin{equation} \label{aux  cent fix point}
{\bf P_o}(u;x) \stackrel{def}{=}  {\bf P} \left[\  B_n |f_n(x) - {\bf E} f_n(x)| > u \ \right]  \le \exp \left[ \  - C_7 \nu(u/V(n))   \ \right];
\end{equation}

\vspace{3mm}

\begin{equation} \label{aux  cent est}
{\bf P}_a^o(u) \stackrel{def}{=} \max_{\vec{k}}  {\bf P} \left[\sup_{x \in L(\vec{k})} B_n (f_n(x) - {\bf E} f_n(x)) > u \ \right].
\end{equation}

\vspace{3mm}

 \ Notice that
\begin{equation} \label{loc max estim}
\max_{x \in Q} ||\zeta_n(x)||G\psi_l(2,\infty) \le \frac{C_6}{n  \ h_n},
\end{equation}

\vspace{3mm}

\begin{equation} \label{Difference estim}
 ||\zeta_n(x) - \zeta_n(y)||G\psi_l(2,\infty) \le \frac{C_8 \ ||x - y|| }{ h_n},
\end{equation}

\vspace{3mm}

and that  when $ \ n \ge 2 \ $

$$
V(n) \stackrel{def}{=} n \ h_n \sim C_7  \ n^{(\beta + d)/(2 \beta + d)} \ (\ln n)^{\beta/(2 \beta + d)}  \to \infty, \ n \to \infty.
$$

\vspace{3mm}

 \ We propose attracting the main result of the section 3.4 of the monograph   \cite{Ostrovsky1999} that

\vspace{3mm}

\begin{equation} \label{Cube estim}
{\bf P}_a^o(u)  \le \exp \left( \ C_4 \ln n - \nu(u \b V(n)   \ \right), \ u \ge u_0 = \const \ge e.
\end{equation}

\vspace{3mm}

 \ Following  we deduce after mild calculations, when $ \ u > e \ $

$$
{\bf P_z}(u) \le \sum_{\vec{k}} \ \exp \left( \ C_5(\beta,d,L) \ln n -\nu (u \ V(n)) \ \right) \le
$$

$$
  n^d \ \exp \left( \ C_6(\beta,d,L) \ \ln n -\nu(u \ V(n)  \ \right) \le
$$

$$
\sup_n \left\{ n^d \ \exp \left( \ C_6(\beta,d,L) \ \ln n -\nu(u \ V(n))  \ \right) \ \right\} \ \le
$$

$$
 C_9(D,\beta,d,L) \exp \left(- \ \nu(u/C_{10}(D,\beta,d,L)   \ ) \ \right), \ C_{9,10}(D,\beta,d,L) \in (0,\infty).
$$

\vspace{3mm}

 \ Ultimately,  the proposition (\ref{main res}) is proved. \\

\vspace{5mm}

{\emph{\textbf{\footnotesize Acknowledgements}.} {\footnotesize M.R. Formica is member of Gruppo
Nazionale per l'Analisi Matematica, la Probabilit\`{a} e le loro Applicazioni (GNAMPA) of the Istituto Nazionale di Alta Matematica (INdAM) and
member of the UMI group \lq\lq Teoria dell'Approssimazione e Applicazioni (T.A.A.)\rq\rq and is partially supported by the INdAM-GNAMPA project,
{\it Risultati di regolarit\`{a} per PDEs in spazi di funzione non-standard}, codice CUP\_E53C22001930001.

\end{document}